\documentclass[12TP,draft]{article}

\begin{document}

\begin{center}
\LARGE\noindent\textbf{ On pre-Hamiltonian Cycles in Hamiltonian Digraphs}\\

\end{center}
\begin{center}
\noindent\textbf{Samvel Kh. Darbinyan}\\

Institute for Informatics and Automation Problems, Armenian National Academy of Sciences

E-mail: samdarbin@ipia.sci.am\\
\end{center}

\textbf{Abstract}

Let $D$ be a strongly connected directed graph of order $n\geq 4$. In \cite{[14]} (J. of Graph Theory, Vol.16, No. 5,  51-59, 1992) Y. Manoussakis proved the following theorem: Suppose that $D$  satisfies the following condition for every triple $x,y,z$ of vertices such that $x$ and $y$ are non-adjacent: If there is no arc from $x$ to $z$, then $d(x)+d(y)+d^+(x)+d^-(z)\geq 3n-2$. If there is no arc from $z$ to $x$, then $d(x)+d(y)+d^-(x)+d^+(z)\geq 3n-2$. Then $D$ is Hamiltonian. In this paper we show that: If $D$ satisfies the condition of Manoussakis' theorem, then $D$ contains a pre-Hamiltonian cycle (i.e., a cycle of length $n-1$) or $n$ is even and $D$ is isomorphic to the complete bipartite digraph with partite sets of cardinalities $n/2$ and $n/2$.\\ 

\textbf{Keywords:} Digraph, cycles, Hamiltonian cycles,  pre-Hamiltonian cycles, longest non-Hamiltonian cycles. \\

\section {Introduction} 

The directed graph (digraph) $D$ is Hamiltonian  if it contains a Hamiltonian cycle,
i.e., a cycle of length $n$ and is pancyclic if it contains cycles of all lengths $m$, $3\leq m\leq n$, where $n$ is the number of vertices in $D$. We recall the following well-known degree conditions (Theorems 1.1-1.8) that guarantee that a digraph is Hamiltonian. In each of the conditions (Theorems 1.1-1.8) below $D$ is a strongly connected digraph of order $n$ :\\

 \textbf{Theorem 1.1} (Ghouila-Houri \cite{[12]}). If $d(x)\geq n$ for all vertices $x\in V(D)$, then $D$ is Hamiltonian.\\

 \textbf{Theorem 1.2} (Woodall \cite{[18]}). If $d^+(x)+d^-(y)\geq n$ for all pairs of vertices $x$ and $y$ such that there is no arc from $x$ to $y$, then $D$ is Hamiltonian.\\
 
\textbf{Theorem 1.3} (Meyniel \cite{[15]}). If $n\geq 2$ and $d(x)+d(y)\geq 2n-1$ for all pairs of non-adjacent
vertices in $D$, then $D$ is Hamiltonian .\\

It is easy to see that Meyniel's theorem is a common generalization of Ghouila-Houri's and Woodall's theorems. For a short proof of Theorem 1.3, see \cite{[5]}. \\

C. Thomassen \cite{[17]} (for $n=2k+1$) and S. Darbinyan \cite{[7]}
(for $n=2k$) proved the following:

\textbf{Theorem 1.4}. If $D$ is a digraph of order
$n\geq 5$ with minimum degree at least $n-1$ and with minimum
semi-degree at least $n/2-1$, then $D$ is Hamiltonian (unless some
extremal cases which are characterized).\\

For the next theorem we need the following:

\textbf{Definition 1} \cite{[14]}. Let $k$ be an arbitrary nonnegative integer. A digraph $D$ satisfies the condition $A_k$ if and only if for every triple $x,y,z$ of vertices such that $x$ and $y$ are non-adjacent: If there is no arc from $x$ to $z$, then $d(x)+d(y)+d^+(x)+d^-(z)\geq 3n-2+k$. If there is no arc from $z$ to $x$, then $d(x)+d(y)+d^-(x)+d^+(z)\geq 3n-2+k$.  \\

\textbf{Theorem 1.5} (Y. Manoussakis \cite{[14]}). If a digraph $D$  satisfies the condition $A_0$, then $D$ is Hamiltonian. \\

Each of these theorems  imposes a degree condition on all pairs of non-adjacent vertices (or on all vertices). In the following three theorems imposes a degree condition only for some pairs of non-adjacent vertices.\\
 
\textbf{Theorem 1.6} (Bang-Jensen, Gutin, H.Li \cite{[2]}). Suppose that $min\{d(x),d(y)\}\geq n-1$ and  $d(x)+d(y)\geq 2n-1$ for any pair of non-adjacent vertices $x,y$ with a common in-neighbour, then $D$ is Hamiltonian.\\

\textbf{Theorem 1.7} (Bang-Jensen, Gutin, H.Li \cite{[2]}). Suppose that $min\{d^+(x)+d^-(y),d^-(x)+d^+(y)\}\geq n$ for any pair of non-adjacent vertices $x,y$ with a common out-neighbour or a common in-neighbour, then $D$ is Hamiltonian.\\

\textbf{Theorem 1.8}  (Bang-Jensen, Guo, Yeo \cite{[3]}). Suppose that $d(x)+d(y)\geq 2n-1$ and  $min\{d^+(x)+d^-(y),d^-(x)+d^+(y)\}\geq n-1$ for any pair of non-adjacent vertices $x,y$ with a common out-neighbour or a common in-neighbour, then $D$ is Hamiltonian.\\

Note that Theorem 1.8 generalizes Theorem 1.7. \\

In \cite{[11], [16], [6], [8]} it was shown that if a digraph $D$ satisfies the condition one of Theorems 1.1, 1.2, 1.3 and 1.4, respectively, then $D$ also is pancyclic (unless some extremal cases which are characterized). It is natural to set the following problem: Characterize those digraphs which satisfy the conditions of Theorem 1.6 (1.7, 1.8) but are not pancyclic. In the many papers (as well as, in the mentioned papers), the existence of a pre-Hamiltonian cycle (i.e., a cycle of length $n-1$) is essential to the show that a given digraph (graph) is pancyclic or not. This indicates that the existence of a pre-Hamiltonian cycle in the a digraph (graph) makes the pancyclic problem  significantly easer, in a sense. In \cite{[9]} the following results were proved: 

(i) if the minimum semi-degree of $D$ at least two and $D$ satisfies the condition of Theorem 1.6 or 

(ii) $D$ is not directed cycle and satisfies the condition of Theorem 1.7, then either $D$ contains a pre-Hamiltonian cycle (i.e., a cycle of length $n-1$) or $n$ is even and $D$ is isomorphic to the complete bipartite digraph or to the complete bipartite digraph minus one arc with partite sets of cardinalities $n/2$ and $n/2$.

 In \cite{[10]} proved that if $D$ is not a directed cycle and satisfies the condition of Theorem 1.8, then $D$ contains a pre-Hamiltonian cycle or a cycle of length  $n-2$.\\

In \cite{[14]} the following conjecture was proposed:

\textbf{Conjecture 1.9}. Any strongly connected digraph satisfying the condition $A_3$ is pancyclic.\\

In this paper using some claims of the proof of Theorem 1.5 (see \cite{[14]}) we prove the following:\\

 \textbf{Theorem 1.10}. Any strongly connected digraph $D$ on $n\geq 4$ vertices satisfying the condition $A_0$ contains a pre-Hamiltonian cycle or $n$ is even and $D$ is isomorphic to the complete bipartite digraph with partite sets of cardinalities of $n/2$ and $n/2$. \\

The following examples show the sharpness of the bound $3n-2$ in the theorem. The digraph consisting of the disjoint union of two complete digraphs with one common vertex and the digraph obtained from a complete bipartite digraph after deleting one arc show that the bound $3n-2$ in the above theorem is best possible. \\

\section {Terminology and Notations}

We shall assume that the reader is familiar with the standard
terminology on the directed graphs (digraph)
 and refer the reader to the monograph of Bang-Jensen and Gutin \cite{[1]} for terminology not discussed here.
  In this paper we consider finite digraphs without loops and multiple arcs. For a digraph $D$, we denote
  by $V(D)$ the vertex set of $D$ and by  $A(D)$ the set of arcs in $D$. The order of $D$ is the number
  of its vertices. Often we will write $D$ instead of $A(D)$ and $V(D)$. The arc of a digraph $D$ directed from
   $x$ to $y$ is denoted by $xy$. For disjoint subsets $A$ and  $B$ of $V(D)$  we define $A(A\rightarrow B)$ \,
   as the set $\{xy\in A(D) / x\in A, y\in B\}$ and $A(A,B)=A(A\rightarrow B)\cup A(B\rightarrow A)$. If $x\in V(D)$
   and $A=\{x\}$ we write $x$ instead of $\{x\}$. If $A$ and $B$ are two disjoint subsets of $V(D)$ such that every
   vertex of $A$ dominates every vertex of $B$, then we say that $A$ dominates $B$, denoted by $A\rightarrow B$. The out-neighborhood of a vertex $x$ is the set $N^+(x)=\{y\in V(D) / xy\in A(D)\}$ and $N^-(x)=\{y\in V(D) / yx\in A(D)\}$ is the in-neighborhood of $x$. Similarly, if $A\subseteq V(D)$, then $N^+(x,A)=\{y\in A / xy\in A(D)\}$ and $N^-(x,A)=\{y\in A / yx\in A(D)\}$. The out-degree of $x$ is $d^+(x)=|N^+(x)|$ and $d^-(x)=|N^-(x)|$ is the in-degree of $x$. Similarly, $d^+(x,A)=|N^+(x,A)|$ and $d^-(x,A)=|N^-(x,A)|$. The degree of the vertex $x$ in $D$ is defined as $d(x)=d^+(x)+d^-(x)$ (similarly, $d(x,A)=d^+(x,A)+d^-(x,A)$). The subdigraph of $D$ induced by a subset $A$ of $V(D)$ is denoted by $\langle A\rangle$. The path (respectively, the cycle) consisting of the distinct vertices $x_1,x_2,\ldots ,x_m$ ( $m\geq 2 $) and the arcs $x_ix_{i+1}$, $i\in [1,m-1]$  (respectively, $x_ix_{i+1}$, $i\in [1,m-1]$, and $x_mx_1$), is denoted by  $x_1x_2\cdots x_m$ (respectively, $x_1x_2\cdots x_mx_1$). We say that $x_1x_2\cdots x_m$ is a path from $x_1$ to $x_m$ or is an $(x_1,x_m)$-path. For a cycle  $C_k:=x_1x_2\cdots x_kx_1$ of length $k$, the subscripts considered modulo $k$, i.e., $x_i=x_s$ for every $s$ and $i$ such that  $i\equiv s\, (\hbox {mod} \,k)$. A cycle that contains the all vertices of $D$ (respectively, the all vertices of $D$ except one) is a Hamiltonian cycle (respectively, is a pre-Hamiltonian cycle). The concept of the pre-Hamiltonian cycle was given in \cite{[13]}. If $P$ is a path containing a subpath from $x$ to $y$ we let $P[x,y]$ denote that subpath. Similarly, if $C$ is a cycle containing vertices $x$ and $y$, $C[x,y]$ denotes the subpath of $C$ from $x$ to $y$. A digraph $D$ is strongly connected (or, just, strong) if there exists a path from $x$ to $y$ and a path from $y$ to $x$ for every pair of distinct vertices $x,y$.
  For an undirected graph $G$, we denote by $G^*$ the symmetric digraph obtained from $G$ by replacing every edge $xy$ with the pair $xy$, $yx$ of arcs.  $K_{p,q}$ denotes the complete bipartite graph  with partite sets of cardinalities $p$ and $q$.  Two distinct vertices $x$ and $y$ are adjacent if $xy\in A(D)$ or $yx\in A(D) $ (or both). For integers $a$ and $b$, $a\leq b$, let $[a,b]$  denote the set of
all integers which are not less than $a$ and are not greater than
$b$. Let $C$ be a non-Hamiltonian cycle in digraph $D$. An $(x,y)$-path $P$ is a $C$-bypass if $|V(P)|\geq 3$, $x\not=y$ and $V(P)\cap V(C)=\{x,y\}$.

\section { Preliminaries }

The following well-known simple Lemmas 3.1-3.4 are the basis of our results
and other theorems on directed cycles and paths in digraphs. They
will be used extensively in the proofs of our results. \\

\textbf{Lemma 3.1} \cite{[11]}. Let $D$ be a digraph of order $n\geq 3$
 containing a
 cycle $C_m$, $m\in [2,n-1] $. Let $x$ be a vertex not contained in this cycle. If $d(x,C_m)\geq m+1$,
 then  $D$ contains a cycle $C_k$ for all  $k\in [2,m+1]$.  \\

The following lemma is a slight modification of a lemma by Bondy and Tomassen \cite{[5]}.\\

\textbf{Lemma 3.2}. Let $D$ be a digraph of order $n\geq 3$
 containing a
 path $P:=x_1x_2\ldots x_m$, $m\in [2,n-1]$ and let $x$ be a vertex not contained in this path.
  If one of the following conditions holds:

 (i) $d(x,P)\geq m+2$;

 (ii) $d(x,P)\geq m+1$ and $xx_1\notin D$ or $x_mx_1\notin D$;

 (iii) $d(x,P)\geq m$, $xx_1\notin D$ and $x_mx\notin D$,
then there is an  $i\in [1,m-1]$ such that
$x_ix,xx_{i+1}\in D$, i.e., $D$ contains a path $x_1x_2\ldots
x_ixx_{i+1}\ldots x_m$ of length $m$ (we say that  $x$ can be
inserted into $P$ or the path $x_1x_2\ldots x_ixx_{i+1}\ldots x_m$
is extended from $P$ with $x$ ). \fbox \\\\

If in Lemma 3.1 and Lemma 3.2 instead of the vertex $x$ consider a path $Q$, then we get the following Lemmas 3.3 and 3.4, respectively. \\

\textbf{Lemma 3.3}. Let $C_k:=x_1x_2\ldots x_kx_1$, $k\geq 2$, be a non-Hamiltonian cycle in a digraph $D$. Moreover, assume that there exists a path $Q:=y_1y_2\ldots y_r$, $r\geq 1$, in $D-C_k$. If $d^-(y_1,C_k)+d^+(y_r,C_k)\geq k+1$, then for all $m\in [r+1,k+r]$ the digraph $D$ contains a cycle $C_m$ of length $m$  with vertex set $V(C_m)\subseteq V(C_k)\cup V(Q)$.  \fbox \\\\

\textbf{Lemma 3.4}. Let $P:=x_1x_2\ldots x_k$, $k\geq 2$, be a non-Hamiltonian path in a digraph $D$. Moreover, assume that there exists a path $Q:=y_1y_2\ldots y_r$, $r\geq 1$, in $D-P$. If $d^-(y_1,P)+d^+(y_r,P)\geq k+d^-(y_1,\{x_k\})+d^+(y_r,\{x_1\})$, then $D$ contains a path from $x_1$ to $x_k$ with vertex set $V(P)\cup V(Q)$. \fbox \\\\

For the proof of our result we also need the following
 
\textbf{Lemma 3.5} \cite{[14]}. Let $D$ be a digraph on $n\geq 3$
vertices satisfying the condition $A_0$. Assume that there are two distinct pairs of non-adjacent vertices $x,y$ and $x,z$  in $D$. 
Then either $d(x)+d(y)\geq 2n-1$ or $d(x)+d(z)\geq 2n-1$.\\

\section {The proof of Theorem 1.10}

In the proof of Theorem 1.10 we often will use the following definition:

\textbf{Definition 2}. Let $P_0:=x_1x_2\ldots x_m$, $m\geq 2$, be an arbitrary $(x_1,x_m)$-path in a digraph $D$ and let 
$y_1,y_2,\ldots y_k\in V(D)- V(P_0)$. For $i\in [1,k]$ we denote by $P_i$ an $(x_1,x_m)$-path in $D$ with vertex set $V(P_{i-1})\cup \{y_j\}$ (if it exists), i.e, $P_i$ is extended path obtained from $P_{i-1}$ with some vertex $y_j$, where $y_j\notin V(P_{i-1}$). If $e+1$ is the maximum possible number of these paths $P_0, P_1,\ldots , P_e$, $e\in [0,k]$, then we say that $P_e$ is extended path obtained from $P_0$ with vertices $y_1,y_2,\ldots , y_k$ as much as possible. Notice that $P_i$ ($i\in [0,e]$) is an $(x_1,x_m)$-path of length $m+i-1$.\\

\textbf{Proof of Theorem 1.10}.
Let $C:=x_1x_2\ldots x_kx_1$ be a longest non-Hamiltonian cycle in $D$ of length $k$, and let $C$ be chosen so that $\langle V(D)- V(C)\rangle$ has the minimum number of connected components. Suppose that $k\leq n-2$ and $n\geq 5$ (the case $n=4$ is trivial). It is easy to show that $k\geq 3$. We will prove that $D$ is isomorphic
to the complete bipartite digraph $K^*_{n/2,n/2}$. Put $R:=V(D)-V(C)$. Let $R_1, R_2,\ldots , R_q$ be the connected components of $\langle R\rangle$ (i.e., if $q\geq 2$, then for any pair $i,j$, $i\not=j$, there is no arc between $R_i$ and $R_j$). In \cite{[14]} it was proved that for any $R_i$, $i\in [1,q]$, the subdigraph $\langle V(C)\cup V(R_i)\rangle$  contains a $C$-bypass. (The existence of a $C$-bypass also follows from Bypass Lemma (see \cite{[4]}), since  $\langle V(C)\cup V(R_i)\rangle$ is strong and condition $A_0$ implies that the underlying graph of the subdigraph $\langle V(C)\cup V(R_i)\rangle$ is 2-connected). Let $P:=x_my_1y_2\ldots y_{t_i}x_{m+\lambda_i}$ be a $C$-bypass in $\langle V(C)\cup V(R_i)\rangle$ ($i\in [1,q]$ is arbitrary) and $\lambda_i$ is considered to be minimum in the sense that there is no $C$-bypass $x_au_1u_2\ldots u_{l_i}x_{a+r_i}$ in $\langle V(C)\cup V(R_i)\rangle$ such that $r_i<\lambda_i$ and $\{x_a,x_{a+r_i}\}$ is a subset of $\{x_{m},x_{m+1},\ldots , x_{m+\lambda_i}\}$. 

We will distinguish two cases, according as there is a $\lambda_i$, $i\in [1,q]$, such that $\lambda_i=1$ or not.

Assume first that $\lambda_i\geq 2$ for all $i\in [1,q]$. For this case one can show that (the proofs as the same as the proofs of Case 1, Lemma 2.3 and Claim 1 in \cite{[14]}) if $\lambda_i\geq 2$, then  
 $t_i=|R_i|=1$, in $\langle V(C)\rangle$ there is an $(x_{m+\lambda_i},x_m)$-path (say, $P'$) of length $k-2$ with vertex set $V(P')=V(C)-\{z_i\}$, where $z_i\in \{x_{m+1},x_{m+2},\ldots , x_{m+\lambda_i-1}\}$ and $d(y_1)+d(z_i)\leq 2n-2$ (note that $y_1$ and $z_i$ are non-adjacent). From  $|R|\geq 2$ and $|R_i|=1$ (for all $i$) it follows that $q\geq 2$. If $u\in R_2$, then $d(u)=d(u,C)\leq k$ (by Lemma 3.1) and $d(z_1,R)=0$ (by minimality of $q$), in particular, the vertices $z_1$ and $u$ are non-adjacent. Therefore $d(z_1)=d(z_1,C)\leq k$ and $d(z_1)+d(u)\leq 2n-2$. This in connection with $d(y_1)+d(z_1)\leq 2n-2$ contradicts Lemma 3.5. 

Assume second that $\lambda_i =1$ for all $i\in [1,q]$. It is clear that $q=1$. Put $t:=t_1$ and $\lambda :=\lambda_1 =1$. Now for this case first we will prove Claims 1-15.

Observe that if $v_1v_2\ldots v_j$ (maybe, $j=1$) is a path in $\langle R\rangle$ and $x_iv_1\in D$, then $v_jx_{i+j}\notin D$ since $C$ is longest non-Hamiltonian cycle in $D$. We shall use this often, without  mentioning  this explicitly. 

From $\lambda=1$ and the maximality of $C$ it follows the following:

\textbf {Claim 1}. $R=\{y_1,y_2,\ldots ,y_t\}$, i.e., $t=n-k\geq 2$ and $y_1y_2\ldots y_t$ is a Hamiltonian path in $\langle R\rangle$, and if $1\leq i<j-1\leq t-1$, then $y_iy_j\notin D$.  \fbox \\\\  

Claim 1 implies that
$$
d^+(y_1,R)=d^-(y_t,R)=1 \quad \hbox{and if}\quad i\in [1,t-1], \,\,\, \hbox{then}\,\,\,d^+(y_i)\leq i; \eqno (1)
$$
$$
d(y_1,R), \, d(y_t,R)\leq n-k \quad \hbox{and if} \quad i\in [2,t-1], \quad \hbox{then} \quad d(y_i,R)\leq n-k+1. \eqno (2)
$$

\textbf {Claim 2}. (i). If $x_iy_1\in D$, then $d^-(x_{i+1},\{y_1,y_2,\ldots ,y_{t-1}\})=0$;

(ii). If $y_tx_{i+1}\in D$, then $d^+(x_{i},\{y_2,y_3,\ldots ,y_{t}\})=d^+(x_{i-1},\{y_1,y_2,\ldots ,y_{t-1}\})=0$;

(iii). $ d(y_j,C)\leq k$ for all $j\in [1,t]$ and in addition, if $x_iy_1$ and $y_tx_{i+1}\in D$, then $d(y_j,C)\leq k-1$ for all $j\in [2,t-1]$ (by Lemma 3.2(iii) and Claim 2(ii)). \fbox \\\\

  \textbf {Claim 3}. Assume that $\langle R\rangle$ is strong. Then there are no two distinct vertices $x_i, x_j$ ($i, j\in [1, k]$) such that  $d^+(x_i, R)\geq 1$, $d^-(x_j, R)\geq 1$, $|C[x_i,x_j]|\geq 3$,
 $d^-(x_{j-1},R)=0$ (respectively, $d^+(x_{i+1},R)=0$), moreover if $|C[x_i,x_j]|\geq 4$, then 
$A(R, C[x_{i+1},x_{j-2}])=\emptyset$ (respectively, $A(R, C[x_{i+2},x_{j-1}])=\emptyset$).

\textbf{Proof}. Suppose that Claim 3 is false. Without loss of generality assume that $x_ky_f, y_gx_l\in D$ ($l\in [2,k-1]$) $d^-(x_{l-1},R)=0$ and if $l\geq 3$, then $A(R, C[x_{1},x_{l-2}])=\emptyset$. The subdigraph $\langle R\rangle$ contains a $(y_f,y_g)$-path (say $P(y_f,y_g)$) since $R$ is strong. We extend the path $P_0:=C[x_l,x_k]$ with the vertices $x_1, x_2, \dots , x_{l-1}$ as much as possible. Then some vertices $z_1,z_2, \ldots , x_d\in \{x_1, x_2, \dots , x_{l-1}\}$, $ d\in [1,l-1]$, are not on the extended path $P_e$ (for otherwise, it is not difficult to see that by Definition 2 there is an $(x_l,x_k)$-path $P_i$, $i\in [0,e]$, which together with the path $P(y_f,y_g)$ and the arcs $x_ky_f, y_gx_l$ forms a non-Hamiltonian cycle longer than $C$). Therefore, by Lemma 3.2(i), for all $s\in  [1,d]$ the following holds 
$$
d(z_s, C)\leq k+d-1. \eqno (3)
$$  
From  $A(R, C[x_{1},x_{l-2}])=\emptyset$ (if $l\geq 3$), $d^-(x_{l-1}, R)=0$ and Lemma 3.2(ii) it follows that 
$$d(y_1, C)\leq k-l+2\quad  \hbox{and} \quad d(y_t, C)\leq k-l+2
$$ 
since neither $y_1$ nor $y_t$ cannot be inserted into $C[x_{l-1},x_k]$. This together with (2) implies that  
$$ 
 d(y_1)\leq n-l+2 \quad \hbox{and} \quad d(y_t)\leq n-l+2. \eqno (4)
$$
If there exists a $z_s$ such that $d(z_s,R)=0$, then by (3) and (4) we obtain that $d(z_s)+d(y_1)\leq 2n-2$ and $d(z_s)+d(y_t)\leq 2n-2$, which contradicts Lemma 3.5. Assume therefore that there is no $z_s$ such that $d(z_s,R)=0$. Then $d=1$, $z_1=x_{l-1}$, $d^+(x_{l-1}, R)\geq 1$,  $d(x_{l-1}, C)\leq k$ (by (3)) and $D$ contains an $(x_l,x_k)$-path with vertex set $V(C)-\{x_{l-1}\}$. From this it follows that $y_f=y_g$, i.e., 
$$
d^+(x_k, R-\{y_f\})=d^-(x_l, R-\{y_f\})=0. \eqno (5) 
$$
 Therefore $D$ contains a cycle $C'$ of length $k$ with vertex set $V(C)\cup \{y_f\}-\{x_{l-1}\}$, and  the vertices $x_{l-1}$, $y_f$ are non-adjacent. From this, (3), $d^-(x_{l-1},R)=0$ and $d(x_{l-1},\{y_f\})=0$ it follows that $d(x_{l-1})\leq n-1$.
 
Assume first that $y_f\not=y_1$. Let $x_{l-1}y_1\in D$. Then $y_f=y_t$ (by Claim 2(i)) and for the triple of vertices $y_t, x_{l-1}, y_1$  condition $A_0$ holds , since $y_1x_{l-1}\notin D$ and $y_t, x_{l-1}$ are non-adjacent. From $d(x_{l-1}, R- \{y_1\})=0$ and (3) it follows that $d(x_{l-1})\leq k+1$. Since $D$ contains no cycle of length $k+1$, it follows that for  the arc $x_{l-1}y_1$ and the cycle $C'$, by Lemma 3.3 the following holds $d^-(x_{l-1},C')+d^+(y_1,C')\leq k$. This together with $d^+(y_1,R)=1$ and $d^-(x_{l-1},R)=0$ implies that $d^-(x_{l-1})+d^+(y_{1})\leq n-2$ (here we consider the cases $k=n-2$ and $k\leq n-3$ separately). Therefore, by condition $A_0$, (4), $d(x_{l-1})\leq n-1$, $l\geq 2$ and $k\leq n-2$, we have
$$
3n-2\leq d(y_t)+d(x_{l-1})+ d^-(x_{l-1})+d^+(y_1)\leq 3n-3,
$$
a contradiction. Let now $x_{l-1}y_1\notin D$. Then the vertices $x_{l-1}$, $y_1$ are non-adjacent and $t\geq 3$ since $d^+(x_{l-1},R)\geq 1$. Using (2) and Lemma 3.2(iii) it is not difficult to see that $d(y_1)\leq n-l$, since $x_ky_1\notin D$ and $y_1x_l\notin D$ (by (5)). Notice that 
$$
d(x_{l-1})=d(x_{l-1},C)+d(x_{l-1},R- \{y_1,y_f\})\leq k+d(x_{l-1},R- \{y_1,y_f\})\leq n-2,
$$
and (by Lemma 3.2(i))
$$ d(y_f)=d(y_f,C)+d(y_f,R)\leq k-l+2+d(y_f,R).
$$
From the last three inequalities we obtain that 
$$
d(y_1)+d(x_{l-1})\leq 2n-2-l 
$$
and
$$
d(y_f)+d(x_{l-1})\leq 2k-l+2+d(x_{l-1},R- \{y_1,y_f\})+d(y_{f},R).
$$
Notice that $$d(x_{l-1},R- \{y_1,y_f\})+d(y_{f},R)\leq n-k-2+n-k=2n-2k-2$$
since if $x_{l-1}y_j\in D$, then $y_jy_f\notin D$, where $y_j\not= y_1,y_f$. Therefore $d(y_f)+d(x_{l-1})\leq 2n-l\leq 2n-2$. This together with $d(y_1)+d(x_{l-1})\leq 2n-2-l$ contradicts Lemma 3.5.

Assume next that $y_f=y_1$. If $x_{l-1}, y_t$ are non-adjacent, then $d(x_{l-1},R)\leq n-k-2$ since $d(x_{l-1},\{y_1,$ $y_t\})$ $=0$ and hence by (3) and $d=1$, $d(x_{l-1})\leq n-2$. Therefore, using (4) we get that $d(y_1)+d(x_{l-1})\leq 2n-2$ and $d(y_t)+d(x_{l-1})\leq 2n-2$ which contradicts Lemma 3.5, since $y_1,x_{l-1}$ and $y_t,x_{l-1}$ are two distinct pairs of non-adjacent vertices. So, we can assume that $x_{l-1}y_t\in D$. Since $C'$ is a longest non-Hamiltonian cycle,  $d^-(x_{l-1},R)=0$ and $d^+(y_t,R-\{y_1\})\leq n-k-2$, from Lemma 3.3 it follows that $d^-(x_{l-1})+d^+(y_t)\leq n-2$. Then from (4) and $d(x_{l-1})\leq n-1$, by condition $A_0$, for the triple of the vertices $x_{l-1}, y_1, y_t$ we obtain that 
$$
3n-2\leq d(y_1)+d(x_{l-1})+d^+(y_t)+d^-(x_{l-1})\leq 3n-l-1\leq 3n-3,
$$
which is a contradiction. Claim 3 is proved. \fbox \\\\

Now we divide the proof of the theorem into two parts: $k\leq n-3$ and $k=n-2$.

Part 1. $k\leq n-3$, i.e., $t\geq 3$. For this part first we will prove the following Claims 4-9 below.\\

\textbf {Claim 4}. Let $t\geq 3$ and $y_ty_1\in D$. Then (i) if $x_iy_1D$, then $d^-(x_{i+2},R)=0$; (ii) if $y_tx_i\in D$, then $d^+(x_{i-2},R)=0$, where $i\in [1,k]$.

\textbf{Proof}. (i). Suppose, on the contrary, that for some $i\in [1,k]$ $x_iy_1\in D$ and $d^-(x_{i+2},R)\not=0$. Without loss of generality, we assume that $x_i=x_1$. Then $d^-(x_3,R-\{y_1\})=0$ and $y_1x_3\in D$. It is easy to see that $y_1$, $x_2$ are non-adjacent and 
$$d^-(x_2, \{y_1,y_2,\ldots , y_{t-1}\})=d^+(x_2, \{y_1,y_3,y_4,\ldots , y_t\})=0, \quad \hbox{i.e.,}\quad d(x_2,R)\leq 2.
\eqno (6)
$$
Since neither $y_1$ nor $x_2$ cannot be inserted into $C[x_3,x_1]$, using (2), (6) and Lemma 3.2, we obtain that 
$$
d(y_1)=d(y_1,C)+d(y_1,R)\leq k+n-k =n \quad \hbox{and}\quad d(x_2)=d(x_2,C)+d(x_2,R)\leq k+2.
$$
On the other hand, by Lemma 3.3 and (1) we have that $d^-(y_t)+d^+(y_1)\leq k+2$ since the arc $y_ty_1$ cannot be inserted into $C$. Therefore, by condition $A_0$, the following holds
$$
3n-2\leq d(y_1)+d(x_2)+d^-(y_t)+d^+(y_1)\leq n+2k+4,
$$
since $y_1, x_2$ are non-adjacent and $y_1y_t\notin D$. 
From  this and $k\leq n-3$ it follows that $k=n-3$, $x_2y_2,y_2y_1\in D$ and hence, the cycle $x_2y_2y_1x_3x_4\ldots x_kx_1x_2$ has length $k+2$, which is a contradiction. 

To show that (ii) is true, it is sufficient to apply the same arguments to the converse digraph of $D$. Claim 4 is proved. \fbox \\\\

\textbf {Claim 5.} If $t\geq 3$ and the vertices $y_1$, $y_t$ are non-adjacent, then $t=3$ and $y_3y_2$,
 $y_2y_1\in D$.

\textbf{Proof.}  Without loss of generality, we assume that $x_1y_1$, $y_tx_2\in D$ (since $\lambda =1$).
 
Assume that 
$y_ty_i\in D$ for some $i\in [2, t-2]$. Then $t\geq 4$. Since neither the arc $y_ty_i$ nor any vertex $y_j$, $j\in [1,t]$,  cannot be inserted into $C$, using Lemma 3.1 and Lemma 3.3, we obtain that 
$$  
d(y_j,C)\leq k \quad \hbox{and} \quad d^-(y_t,C)+d^+(y_i,C)\leq k.    \eqno (7) 
$$ 
 From Claim 1 and the condition that $y_1, y_t$ are non-adjacent it follows that 
$$d(y_1, R)\leq n-k-1 \quad \hbox{and} \quad d(y_t,R)\leq n-k-1. 
$$
From this, since $d(y_j,C)\leq k$ for all $j\in [1,t]$ (by (7)), we obtain that $d(y_1)$ and $d(y_t)\leq n-1$. Now
using (1), (7) and apply condition $A_0$ to the triple of the vertices $y_1, y_t, y_i$, we obtain that
$$ 
3n-2\leq d(y_1)+d(y_t)+d^-(y_t,C)+d^+(y_i,C)+d^-(y_t,R)+d^+(y_i,R)\leq 3n-3,
$$
which is a contradiction. Therefore, if $t\geq 4$, then $y_ty_i\notin D$ for all $i\in [2,t-2]$.

 In a similar way we can also show that $y_iy_1\notin D$ for all $i\in [3, t-1]$. Hence
$$ d(y_1,R)\leq 2, \quad d(y_t,R)\leq 2 \quad \hbox{and} \quad d(y_1)+d(y_t)\leq 2k+4,\eqno (8)
$$ 
since $d(y_i)\leq k$ for all $i\in [1,t]$.

If $t\geq 4$, then $y_1,y_t$ and $y_1, y_{t-1}$ are two distinct pairs of non-adjacent vertices. From (8) and $k\leq n-4$ it follows that $d(y_1)+d(y_t)\leq 2n-4$. On the other hand, since $d(y_1)\leq k+2$, $d(y_{t-1},C)\leq k-1$ (by Claim 2 and Lemma 3.2(iii)) and  $d(y_{t-1},R)\leq n-k$ (by Claim 1), we have that
$$
d(y_1)+d(y_{t-1})\leq 2n-3.
$$
This together with  $d(y_1)+d(y_t)\leq 2n-4$ contradicts Lemma 3.5. Therefore $t=3$.

Now we show that $y_3y_2\in D$. Assume that this is false, i.e., $y_3y_2\notin D$. Then we can apply condition $A_0$ to the triple of the vertices $y_1,y_3,y_2$, since the vertices $y_1,y_3$ are non-adjacent and $y_3y_2\notin D$. Notice that the arc $y_2y_3$ cannot be inserted into $C$ and hence $d^-(y_2,C)+d^+(y_3,C)\leq k$ (by Lemma 3.3). Therefore by $A_0$ and  Claim 2, we obtain that
$$ 
3n-2\leq d(y_1)+d(y_3)+d^-(y_2)+d^+(y_3)\leq 3k+4\leq 3n-5,
$$
which is a contradiction. Therefore $y_3y_2\in D$.

In a similar way, as above, we can show that $y_2y_1\in D$. Claim 5 is proved. \fbox \\\\

\textbf {Claim 6.} If $t\geq 3$, then  $y_ty_1\in D$.
 
\textbf{Proof.} Suppose, on the contrary, that $t\geq 3$ and $y_ty_1\notin D$, i.e., $y_1,y_t$ are non-adjacent. Then by Claim 5, $t=3$ and $y_3y_2, y_2y_1\in D$. Without loss of generality, assume that $x_1y_1$ and $y_3x_2\in D$ (since $\lambda =1$). Notice that $d(y_1),d(y_3)\leq n-1$ (by Lemma 3.1). We will distinguish two cases, according as there is an arc from $R$ to $\{x_3,x_4,\ldots , x_k\}$ or not.

\textbf{Case 6.1.} $A(R\rightarrow \{x_3,x_4,\ldots , x_k\})\not= \emptyset$. Then there exists a vertex $x_l$ with $l\in [3,k]$ such that $d^-(x_l, R)\geq 1$ and $A(R\rightarrow \{x_3,x_4,\ldots , x_{l-1}\})= \emptyset$. 

If $l=3$, then from $d^-(x_3,\{y_2,y_3\})=0$ it follows that $y_1x_3\in D$. From this it is easy to see that $d(x_2,\{y_1,y_2\})=0$. Since neither $y_1$ nor $y_3$ and nor $x_2$ cannot be inserted into $C[x_3,x_1]$ using Lemma 3.2 we obtain that $d(y_1)$, $d(y_3)$ and $d(x_2)\leq n-1$. Hence, $d(y_1)+d(y_3)\leq 2n-2$ and $d(y_1)+d(x_2)\leq 2n-2$,
which contradicts Lemma 3.5. 

Assume therefore that $l\geq 4$. From Claim 3, $x_1y_1\in D$ and the minimality of $l$ it follows that $d^+(x_{l-1},R)\geq 1$. Without loss of generality, we may assume that $y_gx_{l}\in D$ and $x_{l-1}y_f$. It is easy to see that $y_f\not=y_g$, $y_f,y_g\in \{y_1,y_3\}$ and the vertices $x_{l-1}, x_g$ are non-adjacent. 

Assume first that $l=4$. Then it is easy to see that $y_g=y_1$ and $y_f=y_3$, i.e., $y_1x_4$ and $x_3y_3\in D$. Then clearly the vertices  $x_2,y_2$ are non-adjacent and $x_2y_3\notin D$. Therefore $x_2y_1\notin D$ (for otherwise if $x_2y_1\in D$, then Claim 3 is not true since $d^-(x_3,R)=0$). Therefore $d(x_2,\{y_1,y_2\})=0$. Notice that $x_2$
cannot be inserted into the path $C[x_4,x_1]$ (for otherwise in $D$ there is a cycle of length $n-3$ for which Claim 5 is not true since $y_3x_3\notin D$). Now by Lemma 3.2 and the above observation we obtain that 
$$d(x_2)=d(x_2,C[x_4,x_1])+d(x_2,R)+d(x_2,\{x_3\})\leq n-1.$$
 Therefore 
 $d(y_1)+d(x_2)\leq 2n-2$, which together with $d(y_1)+d(y_3)\leq 2n-2$
 contradicts Lemma 3.5, since $y_1,x_2$ and $y_1,y_3$ are two distinct pairs of non-adjacent vertices. 

Assume next that $l\geq 5$. From $x_1y_1\in D$, $d^-(x_{l-1}, R)=0$ and Claim 3 it follows that 
$A( \{x_2,x_3,\ldots ,$ $ x_{l-2}\}\rightarrow R)=\emptyset$, in particular, $x_2y_3\notin D$. Therefore 
 $A( \{x_3,x_4,\ldots , x_{l-2}\},R)=\emptyset$, $d(x_2,R)= 1$ (only $y_3x_2\in D$), $d(x_{l-1},R)= 1$ and
$$
d(y_1,\{x_2,x_3,\ldots , x_{l-2}\})=d(y_3,\{x_3,x_4,\ldots , x_{l-2}\})=0. \eqno (9) 
$$
Since neither $y_1$ nor $y_3$ cannot be inserted into $C$, $x_2y_3\notin D$ and $d^-(x_{l-1}, R)=0$, using (9) and Lemma 3.2 we obtain that
 $d(y_1)$ and $d(y_3)\leq k-l+5$. Therefore $d(y_1)+d(y_3)\leq 2n-6$. Now we extend the path  $P_0:=C[x_l,x_1]$ with the vertices 
$x_2,x_3,\ldots , x_{l-1}$ as much as possible. Then some vertices $z_1,z_2,\ldots , z_d\in \{x_2,x_3,\ldots , x_{l-1}\}$, $d\in [1,l-2]$, are not  on the extended path $P_e$. Therefore $d(z_i,C)\leq k+d-1$ and hence, $d(z_i)\leq k+d$ for all $i\in [1,d]$. It is not difficult to show that there is a $z_i$ which is not adjacent with $y_1$. Thus we have $d(y_1)+d(z_i)\leq 2n-3$. This together with $d(y_1)+d(y_3)\leq 2n-6$ contradicts Lemma 3.5 since $y_1,z_i$ and $y_1,y_3$ are  two distinct pairs of non-adjacent vertices. In each case we have a contradiction, and hence the discussion of Case 6.1 is completed.

\textbf{Case 6.2.} $A(R\rightarrow \{x_3,x_4,\ldots ,$  $ x_k\})= \emptyset$. Without loss of generality, we may assume that 
$A( \{x_3,x_4,$ $\ldots , x_k\}\rightarrow R)= \emptyset$ (for otherwise, we consider the converse digraph of $D$ for which the considered Case 6.1 holds). Therefore $A(R,\{x_3,x_4,\ldots , x_k\})=\emptyset$. In particular, $x_k$ is not adjacent with the vertices $y_1$ and $y_3$. Notice that 
$$
d(y_1)=d(y_1,R)+d(y_1,C)\leq 2+d(y_1,\{x_1,x_2\})\leq 5,
$$ 
$d(y_3)\leq 5$ and $d(x_k)=d(x_k,C)\leq 2n-8$. Therefore 
$d(x_k)+d(y_1)\leq 2n-3$ and $d(x_k)+d(y_3)\leq 2n-3$, which contradicts Lemma 3.5. Claim 6 is proved. \fbox \\\\

\textbf {Claim 7.} If $t\geq 3$ and for some $i\in [1,k]$ $x_iy_1$, then $A(R\rightarrow C[x_{i+2},x_{i-1}])=\emptyset$. 

\textbf{Proof}. Suppose that the claim is not true. Without loss of generality,  we may assume that $x_1y_1\in D$ and  $A(R\rightarrow \{x_3,x_4,\ldots ,x_k\})\not=\emptyset$. Then there is a vertex $x_l$ with $l\in [3,k]$ such that $d^-(x_l, R)\geq 1$ and if $l\geq 4$, then $A(R\rightarrow \{x_3,x_4,\ldots ,x_{l-1}\})=\emptyset$. We have that $y_ty_1\in D$ (by Claim 6). In particular, $y_ty_1\in D$ implies that  $\langle R \rangle$ is strong. On the other hand, by Claim 4(i),  $d^-(x_3,R)=0$ and hence, $l\geq 4$. From $x_1y_1\in D$ it follows that there exists a vertex $x_r$ with $r\in [1,l-1]$ such that $d^+ (x_r,R)\geq 1$. Choose $r$ with these properties as maximal as possible. Let $x_ry_f$ and $y_gx_l\in D$. Notice that in $\langle R\rangle$ there is a $(y_f,y_g)$-path since $\langle R\rangle$ is strong. Using  Claim 3 we obtain that $r=l-1$. Then $y_f\not=y_g$ and in $\langle R \rangle$ any $(y_f,y_g)$-path is a Hamiltonian path. Since $\langle R\rangle$ is strong, from $d^-(x_{l-1}, R)=0$, $d^-(x_l,R)\geq 1$ and from Claim 3 it follows that $A(\{x_{2},x_{3},\ldots , x_{l-2}\}\rightarrow R)= \emptyset$, in particular, $d^+(x_2,R)=0$. Then
$$
A(\{x_{3},x_{4},\ldots , x_{l-2}\}, R)= \emptyset, \quad  d(y_1, \{x_2,x_3,\ldots ,x_{l-2}\})=d(x_2, \{y_1,y_2,\ldots ,y_{t-1}\})=0. \eqno (10)
$$
Note that $x_2$, $y_1$ and $x_2$, $y_2$ are two distinct pairs of non-adjacent vertices. We extend the path $P_0:=C[x_l,x_1]$ with the vertices $x_2,x_3,\ldots , x_{l-1}$ as much as possible. Then some vertices $z_1,z_2, \ldots ,z_d \in \{x_2,x_3,\ldots , x_{l-1}\}$, 
where $d\in [1,l-2]$, are not on the extended path $P_e$ (for otherwise, since in $\langle R\rangle$ there is a $(y_f,y_g)$-path, using the path $P_{e-1}$ or $P_e$ we obtain a non-Hamiltonian cycle longer than $C$). By Lemma 3.2, for all $i\in [1,d]$ we have that 
$$
d(z_i,C)\leq k+d-1 \quad \hbox{and} \quad d(z_i)=d(z_i,C)+d(z_i,R)\leq k+d-1+d(z_i,R). \eqno (11)
$$
Assume that there is a vertex $z_i\not= x_{l-1}$. Then, by (10), $d(z_i,R)\leq 1$ (since $d(x_2,R)\leq 1$). Notice that $y_1$, $z_i$ and $y_2$, $z_i$ are two distinct pairs of non-adjacent vertices (by (10)).  Since neither $y_1$ nor $y_2$ cannot be inserted into $C[x_{l-1},x_1]$ and $y_1x_{l-1}\notin D$, $y_2x_{l-1}\notin D$, by Lemma 3.2(ii) and (10) for  $j=1$ and $2$ we obtain that 
$$
d(y_j,C)=d(y_j,C[x_{l-1},x_1])\leq k-l+3. \eqno (12)
$$
In particular, by (2), 
$$d(y_1)=d(y_1,C)+d(y_1,R)\leq k-l+3+n-k= n-l+3.$$
This together with (11) and $d(z_i,R)\leq 1$ implies that 
$$d(y_1)+d(z_i)\leq 2n-2, 
$$
since $k\leq n-3$ and $d\leq l-2$. Therefore, by Lemma 3.5, $d(y_2)+d(z_i)\geq 2n-1$. Hence, by (2) and (11) we have 
$$ 2n-1\leq d(y_2)+d(z_i)\leq n+d +d(z_i,R)+d(y_2,C).$$ 
From this and (12) it follows that $d(y_2,C)=k-l+3$, $d(z_i,R)=1$ and $k=n-3$. Then $z_i=x_2$ and $y_tx_2\in D$ (by (10) and $d^+(x_2,R)=0$). Therefore $x_1y_2\notin D$. From this, $y_2x_{l-1}\notin D$ and $d(y_2,C)=k-l+3$, by Lemma 3.2(iii) we conclude that $y_2$ can be inserted into $C$, which contrary to our assumption.

Now assume that there is no $z_i\not=x_{l-1}$. Then $d=1$, $z_1=x_{l-1}$ and $d^-(x_l, \{y_2,y_3,\ldots , y_t\})=0$ (since $x_1y_1\in D$). Therefore $y_1x_l\in D$ and hence, $d(x_{l-1},R-\{y_2\})=0$ (since $y_ty_1\in D$ and $l$ is minimal), in particular, the vertices $y_t, x_{l-1}$ are non-adjacent. This together with (11) implies that $d(x_{l-1})\leq k+1$ (only $x_{l-1}y_2\in D$ is possible). Notice that neither $y_t$ nor the arc $y_ty_1$ cannot be inserted into $C$, and therefore, by Lemmas 3.2, 3.3 and by (1), (2) we obtain that $d(y_t)\leq n$ and $d^-(y_t)+d^+(y_1)\leq k+2$. Now for the triple of the vertices $y_t, x_{l-1}$, $y_1$, by condition $A_0$, we obtain that 
$$  
3n-2\leq d(x_{l-1})+d(y_t)+d^-(y_t)+d^+(y_1)\leq 3n-3
$$
since $k\leq n-3$, which is a contradiction. Claim 7 is proved. \fbox \\\\

\textbf {Claim 8.} If $t\geq 3$, $x_1y_1$ and $y_tx_2\in D$,  
then $d^-(x_1, R)=0$.

\textbf{Proof}. Assume that $d^-(x_1, R)\geq 1$. By Claim 6, $y_ty_1\in D$. Now using Claims 4(ii) and 7, we obtain that $d^+(x_k,R)=0$ and  
$$A(R\rightarrow \{x_3,x_4,\ldots , x_k\})=\emptyset. \eqno (13) $$
In particular, $d(x_k,R)=0$.
 This together with $d^-(x_1, R)\geq 1$, (13) and Claim 3 implies that $A(\{x_2,x_3,$ $\ldots , x_{k-1}\}\rightarrow R)=\emptyset$. Now again using (13) we get that $A(\{x_3,x_4,\ldots , x_{k}\},R)=\emptyset$. This together with $d^+(x_2,R)=d^-(x_2,\{y_1,y_2,\ldots ,
y_{t-1}\})=0$ implies that $d(x_2,R)=1$, $d(y_2,C)\leq 1$ (only $y_2x_1\in D$ is possible) and $d(x_3,R)=0$. Therefore, by (2),
$$
d(y_2)+d(x_3)=d(y_2,C)+d(y_2,R)+d(x_3,R)+d(x_3,C)\leq n+k\leq 2n-3
$$
and $d(y_2)+d(x_2)\leq 2n-2$, which contradicts Lemma 3.5 since $y_2,x_3$ and $y_2,x_2$ are two distinct pairs of non-adjacent vertices. This completes the proof of Claim 8. \fbox \\\\

\textbf{Claim 9.} If $t\geq 3$, $x_1y_1$ and $y_tx_2\in D$, then $A(\{x_3,x_4,\ldots , x_{k}\}\rightarrow R)=\emptyset$.

\textbf{Proof.} By Claim 6, $y_ty_1\in D$. Suppose that $A(\{x_3,x_4,\ldots , x_{k}\}\rightarrow R)\not=\emptyset$. Recall that Claim 4(ii) implies that
$d^+(x_k,R)=0$.
 Let $x_r$, $r\in [3,k-1]$, be chosen so that $x_ry_i\in D$ for some $i\in [1,t]$ and $r$ is maximum possible. Then $A(\{x_{r+1},x_{r+2},\ldots , x_{k}\}, R)=\emptyset$ and $d^-(x_1,R)=0$ by Claims 7 and 8, respectively. This together with $y_tx_2\in D$ contradicts  Claim 3. Claim 9 is proved. \fbox \\\\

We are now ready to  complete the proof of Theorem 1.10 for Part 1 (when $k\leq n-3$, i.e., $t\geq 3$). By Claim 6, if $t\geq 3$, then $y_ty_1\in D$. Without loss of generality, we may assume that $x_1y_1$ and $y_tx_2\in D$ since $\lambda =1$. Then
from Claims 7, 8 and 9 it follows that
$$
A(R\rightarrow \{x_3,x_4,\ldots , x_k, x_1\})=A(\{x_3,x_4,\ldots , x_{k}\}\rightarrow R)=\emptyset.
$$
From this and
$$
d^-(x_2,\{y_1,y_2,\ldots , y_{t-1}\})=d^+(x_1,\{y_2,y_3,\ldots , y_{t}\})=0
$$
 we obtain that  $x_1,y_2$ and $x_1,y_t$ are two distinct pairs of non-adjacent vertices and $d(y_2,C)\leq 1$, $d(y_t,C)\leq 2$, $d(x_1,R)=1$. Therefore $d(y_2)\leq n-k+2$, $d(y_t)\leq n-k+2$ (by (2))
and $d(x_1)\leq 2k-1$. These inequalities imply that $d(y_2)+d(x_1)\leq 2n-2$ and $d(y_t)+d(x_1)\leq 2n-2$, which contradicts Lemma 3.5. and completes the discussion of Part 1.\\

Part 2. $k=n-2$, i.e., $t=2$. For this part first we will prove Claims 10-15 below.

\textbf{Claim 10}. If  $x_iy_f\in D$ and $y_2y_1\notin D$, where $i\in [1,n-2]$ and $f\in [1,2]$, then there is no $l\in [3,n-2]$ such that $y_fx_{i+l-1}\in D$ and $d(y_f,\{x_{i+1},x_{i+2},\ldots ,x_{i+l-2}\})=0$.

 \textbf{Proof}. The proof is by contradiction. Suppose that $x_iy_f, y_fx_{i+l-1}\in D$ and 
$d(y_f,\{x_{i+1},x_{i+2},\ldots ,$ $x_{i+l-2}\})$ $=0$ for some $l\in [3,n-2]$. Without loss of generality, we may assume that $x_i=x_1$. Then $x_1y_f$, $y_fx_l\in D$ and $d(y_f,\{x_{2},x_{3},\ldots ,$ $x_{l-1}\})$ $=0$. Since $D$ contains no cycle of length $n-1$, using Lemmas 3.2 and 3.3, we obtain that 
$$
d^-(y_1)+d^+(y_2)\leq n-2 \quad\hbox{and} \quad d(y_f)\leq n-l+2. \eqno (14) 
$$
We extend the path $P_0:=C[x_l,x_1]$ with the vertices $x_2, x_3, \ldots , x_{l-1}$ as much as possible. Then some vertices $z_1, z_2, \ldots , z_d\in \{x_2, x_3, \ldots , x_{l-1}\}$, $d\in [1,l-2]$, are not on the extended path $P_e$. Therefore by Lemma 3.2, $d(z_1)=d(z_1,C)+d(z_1,\{y_{3-f}\})\leq n+d-1$. Now, since the vertices $y_f, z_1$ are non-adjacent and $y_2y_1\notin D$, by condition $A_0$ and (14) we have 
$$
3n-2\leq d(y_f)+d(z_1)+d^-(y_1)+d^+(y_2)\leq 3n-3,
$$
a contradiction. Claim 10 is proved. \fbox \\\\

\textbf{Claim 11}. $y_2y_1\in D$ (i.e., if $k=n-2$, then $\langle V(D)-V(C)\rangle$ is strong).

\textbf{Proof}. Suppose, on the contrary, that $y_2y_1\notin D$.
Without loss of generality, we may assume that $x_1y_1\in D$ and the vertices $y_1, x_2$ are non-adjacent. Then $y_2x_3\notin D$ and since $D$ contains no cycle of length $n-1$, using Lemma 3.3 for the arc $y_1y_2$ we obtain that
$$
d^-(y_1)+d^+(y_2)\leq n-2. \eqno (15)
$$ 

\textbf{Case 11.1}. $d^+(y_1,C[x_3,x_{n-2}])\geq 1$. Let $x_l$, $l\in [3, n-2]$, be chosen so that $y_1x_l\in D$ and $l$ is minimum, i.e., $d^+(y_1, C[x_2,x_{l-1}])=0$. It is easy to see that the vertices $y_1$ and $x_{l-1}$ are non-adjacent. By Claim 10 we can assume that $l\geq 5$ (if $l\leq 4$, then $d(y_1,C[x_2,x_{l-1}])=0$, a contradiction to Claim 10) and $d^-(y_1, C[x_3,x_{l-2}])\geq 1$. It follows that there exists a vertex $x_r$ with $r\in [3,l-2]$
such that $x_ry_1\in D$ and $d(y_1,C[x_{r+1},x_{l-1}])=0$. Consequently, for the vertices $y_1$, $x_r$ and $x_l$ Claim 10 is not true, a contradiction.

\textbf{Case 11.2}. $d^+(y_1,C[x_3,x_{n-2}])=0$. Then  $d^+(y_1,C[x_2,x_{n-2}])=0$ and either $y_1x_1\in D$ or $y_1x_1\notin D$.

\textbf{Subcase 11.2.1}. $y_1x_1\in D$. Then $x_{n-2}y_1\notin D$ and hence, the vertices $y_1, x_{n-2}$ are non-adjacent. 
Claim 10 implies that $d^-(y_1, C[x_2,x_{n-2}])=0$. This together with $d^+(y_1, C[x_2,x_{n-2}])=0$ and $y_2y_1\notin D$ gives $d(y_1)=3$. Clearly, $d(x_2)\leq 2n-4$ and hence, for the vertices $y_1, y_2, x_2$ by condition $A_0$ and (15) we have, 
$$
3n-2\leq d(y_1)+d(x_2)+d^-(y_1)+d^+(y_2)\leq 3n-3,
$$
which is a contradiction.

\textbf{Subcase 11.2.2}. $y_1x_1\notin D$. Then $d^+(y_1,C)=0$, $d^+(y_1)=1$ and $d^+(y_2,C)\geq 1$ since $D$ is strong. Without loss of generality, we may assume that $d^-(y_2,C)=0$ (for otherwise for the vertex $y_2$ in the converse digraph of $D$ we would have the above considered Case 11.1 or Subcase 11.2.1). Using Lemma 3.5, it is not difficult to show that $n\geq 6$.

Suppose first that $y_2x_2\in D$. Then $x_{n-2}y_1\notin D$ and hence, the vertices $x_{n-2}, y_1$ are non-adjacent. 

Let for some $l\in [3,n-3]$ $x_ly_1\in D$ and 
$d^-(y_1,C[x_{l+1},x_{n-2}])=0$. Then $d(y_1,C[x_{l+1},x_{n-2}])=0$ and $d(y_1)\leq l$ since $d^+(y_1,C)=0$ and $x_2,y_1$ are non-adjacent. Extend the path $P_0:=C[x_2,x_l]$ with the vertices $x_{l+1}, x_{l+2}, \ldots , x_{n-2}, x_1$ as much as possible. Then some vertices $z_1, z_2, \ldots , z_d\in \{x_{l+1}, x_{l+2}, \ldots ,$ $ x_{n-2}, x_1\}$, $d\in [2,n-l-1]$, are not on the extended path $P_e$. For a vertex $z_i\not=x_1$ by Lemma 3.2 we obtain that $d(z_i)=d(z_i,C)+d(z_i,\{y_2\})\leq n+d-1$. Therefore, since $y_2y_1\notin D$ and the vertices $z_i, y_1$ are non-adjacent, by condition $A_0$ and (15), we get that 
$$
3n-2\leq d(y_1)+d(z_i)+d^-(y_1)+d^+(y_2) \leq 3n-4,
$$
which is a contradiction. 

Let now $x_ly_1 \notin D$ for all $l\in [3,n-2]$, i.e., $d^-(y_1, C[x_3,x_{n-2}])=0$. Then from $d^+(y_1, C[x_2,x_{n-2}])=0$ and $x_{n-2}y_2\notin D$ (since $d^-(y_2,C)=0$) it follows that $d(y_1)=2$ and $d(x_{n-2})\leq 2n-5$. From this, since the vertices $y_1$, $x_{n-2}$ are non-adjacent and $y_2y_1\notin D$, by    condition $A_0$  and (15) we have that
$$
3n-2\leq d(y_1)+d(x_{n-2})+d^-(y_1)+d^+(y_2) \leq 3n-5,
$$        
which is a contradiction.\\

Suppose next that $y_2x_2\notin D$. Then $d(y_2,\{x_2,x_3\})=0$, since $d^-(y_2, C)=0$. Let for some $l\in [4,n-2]$ $y_2x_l\in D$
 and
$d^+(y_2,C[x_2,x_{l-1}])=0$. Then $d(y_2,C[x_2,x_{l-1}])=0$ and the vertices $y_1$, $x_{l-2}$ are non-adjacent since $d^+(y_1,C[x_2,x_{n-2}])=0$. It is easy to see that there exists a vertex $x_r\in \{x_1, x_2,\ldots , x_{l-3}\}$ such that $x_ry_1\in D$ and 
$d(y_1,C[x_{r+1},x_{l-2}])=0$. Thus we have that $A(R,C[x_{r+1},x_{l-2}])=\emptyset$. Notice that $d(y_2)\leq n-l+1$ since $d^-(y_2,C)=0$ and $d(y_2,C[x_2,x_{l-1}])=0$. 
 We extend the path $P_0:=C[x_l,x_r]$ with the vertices $x_{r+1}, x_{r+2}, \ldots , x_{l-1}$ as much as possible. Then some vertices $z_1, z_2, \ldots , z_d\in \{x_{r+1}, x_{r+2}, \ldots , x_{l-1}\}$, $d\in [2,l-r-1]$, are not on the extended path $P_e$. Therefore by Lemma 3.2 for $z_i\not=x_{l-1}$ we have,  $d(z_i)\leq n+d-3$. Now by condition $A_0$ and (15) we obtain that 
$$
3n-2\leq d(y_2)+d(z_{i})+d^-(y_1)+d^+(y_2)<  3n-3,
$$
a contradiction. Let now $d^+(y_2,\{x_2,x_3,\ldots , x_{n-2}\})=0$. Then $d(y_2)=2$, $d(x_2)\leq 2n-6$ and the vertices $x_2,y_2$ are non-adjacent. By condition $A_0$ we have 
$$ 
3n-2\leq d(y_2)+d(x_{2})+d^-(y_1)+d^+(y_2) < 3n-3,
$$
a contradiction. Claim 11 is proved. \fbox \\\\

\textbf{Claim 12}. For any $i\in [1,n-2]$ and $f\in [1,2]$ the following holds

 i) $d^-(y_f,\{x_{i-1},x_i\})\leq 1$ and ii) $d^+(y_f,\{x_{i-1},x_i\})\leq 1$.

\textbf{Proof}. Suppose that the claim is not true. Without loss of generality, we may assume that $x_{n-3}y_1$, $x_{n-2}y_1\in D$ and $y_1,x_1$ are non-adjacent. By Claim 11, $y_2y_1\in D$. It is easy to see that $d^+(y_2, \{x_1, x_2\})=0$, $y_1x_{n-2}\notin D$ and  
 $y_1x_2\notin D$ (for otherwise, if $y_1x_2\in D$, then $x_{n-2}y_1x_2x_3\ldots x_{n-3}x_{n-2}$ is a cycle of length $n-2$ for which $\langle\{y_2,x_1\}\rangle$ is not strong, a contradiction to Claim 11). Therefore, 
$A(R\rightarrow \{x_1,x_2\})=\emptyset$. It is not difficult to check that $n\geq 6$.

Assume first that $A(R\rightarrow \{x_3,x_4,\ldots , x_{n-3}\})\not= \emptyset$. Now let $x_l$,  $l\in [3,n-3]$, be the first vertex after $x_2$ that $d^-(x_l, R)\geq 1$. Then $A(R\rightarrow \{x_1,x_2,\ldots , x_{l-1}\})= \emptyset$ since 
$A(R\rightarrow \{x_1,x_2\})=\emptyset$ (in particular, $d^-(x_{l-1},R)=\emptyset$). From the minimality of $l$ and $x_{n-2}y_1\in D$ it follows that there is a vertex $x_r\in \{x_{n-2},x_1,x_2, \ldots , x_{l-2}\}$ such that $d^+(x_r, R)\geq 1$ and $A(\{x_{r+1},x_{r+2}, \ldots , x_{l-2}\}, R)=\emptyset$ (if $x_r=x_{n-2}$, then $x_{r+1}=x_1$).  This is contrary to Claim 3 since $d^-(x_{l-1},R)=0$ and $\langle R\rangle$ is strong. 
 
Assume next that $A(R\rightarrow \{x_3,x_4,\ldots , x_{n-3}\})= \emptyset$. This together with 
$A(R\rightarrow \{x_1,x_2\})=\emptyset$ gives that $A(R\rightarrow \{x_1,x_2,\ldots ,x_{n-3}\})=\emptyset$. From this,
since $D$ is strong and $y_1x_{n-2}\notin D$, it follows that $y_2x_{n-2}\in D$. Then $x_{n-3}y_2\notin D$ and $x_{n-4}y_1\notin D$. Now using Claim 11 we obtain that $d(y_2,\{x_{n-4},x_{n-3}\})=0$ and $d(x_{n-4},R)=0$. 
If $A(\{x_1, x_2,\ldots , x_{n-5}\}\rightarrow R)\not=\emptyset$, then there is a vertex  $x_r$ with $r\in [1,n-5]$ such that $d^+(x_r, R)\geq 1$ and $A(R,\{x_{r+1},x_{r+2},\ldots , x_{n-4}\})= \emptyset$ ($n\geq 6$) which contradicts Claim 3,  since $y_2x_{n-2}\in D$ and $d^-(x_{n-3}, R)=0$.  
Assume therefore that $A(\{x_1, x_2,\ldots , x_{n-4}\}\rightarrow R)=\emptyset$. Thus we have that $A(\{x_1, x_2,$ $\ldots , x_{n-4}\}, R)=\emptyset$ and $d^-(x_{n-3},R)=0$. Then $d(y_1)= 4$, $d(y_2)\leq 4$ and $d(x_1)\leq 2n-6$. From this it follows that $d(y_1)+d(x_1)\leq 2n-2$ and $d(y_2)+d(x_1)\leq 2n-2$ which contradicts Lemma 3.5. This contradiction proves that $d^-(y_f,\{x_{i-1},x_i\})\leq 1$ for all $i\in [1,n-2]$ and $f\in [1,2]$. Similarly, one can show that $d^+(y_f,\{x_{i-1},x_i\})\leq 1$. Claim 12 is proved. \fbox \\\\

\textbf{Claim 13}. If $x_{i}y_f\in D$ (respectively, $y_fx_i\in D$), then $d(y_f,\{ x_{i+2}\})\not=0$ (respectively, $d(y_f,\{x_{i-2}\})\not=0$), where $i\in [1,n-2]$ and $f\in [1,2]$. 

\textbf{Proof}. Suppose that the claim is not true. By Claim 11, $y_2y_1\in D$. Without loss of generality, we may assume that $x_{n-2}y_1\in D$ and $d(y_1,\{x_2\})=0$, i.e., the vertices $y_1$ and 
$x_2$ are non-adjacent. Claim 12 implies that the vertices $y_1, x_1$ also are non-adjacent. Note that $y_2x_2\notin D$ and hence $d^-(x_2,R)=0$. Now it is not difficult to see that if $n=5$, then  $d(y_1)+d(x_1)\leq 8$ and $d(y_1)+d(x_2)\leq 8$, a contradiction to Lemma 3.5. Assume therefore that $n\geq 6$ and consider the following cases.\\

\textbf{Case 13.1}. $A(R\rightarrow \{x_3,x_4,\ldots , x_{n-3}\})\not= \emptyset$. Then there is a vertex $x_l$ with $l\in [3, n-3]$ such that $d^-(x_l, R)\geq 1$ and $A(R\rightarrow \{x_2,x_3,\ldots , x_{l-1}\})= \emptyset$. We now consider the case $l=3$ and the case $l\geq 4$ separately.

 Assume that $l=3$. Then $y_2x_3\in D$ or $y_1x_3\in D$. 

Let $y_2x_3\in D$. Then the vertices $y_2, x_2$ are non-adjacent. Since the vertices $y_1,x_2$ are non-adjacent Claim 11 implies that $x_1y_2\notin D$. This contradicts Claim 3 because of $d(x_2,R)=0$ and $d^+(x_1,R)=0$. 

Let now $y_1x_3\in D$ and $y_2x_3\notin D$. Then it is easy to see that $x_1y_2\notin D$ and $y_2x_2\notin D$. From this and Claim 11 implies that neither $x_1$ nor $x_2$ cannot be inserted into $C[x_3,x_{n-2}]$. Notice that if $x_2y_2\in D$, then $x_{n-2}x_2\notin D$, and if $y_2x_1\in D$, then $x_1x_3\notin D$. Now using Lemma 3.2, we obtain that $d(y_1)$, $d(x_1)$ and 
$d(x_2)\leq n-1$ since $d(y_1,\{x_1,x_2\})=0$. Therefore $d(y_1)+d(x_1)\leq 2n-2$ and $d(y_1)+d(x_2)\leq 2n-2$, which contradicts Lemma 3.5 since $y_1,x_1$ and $y_1,x_2$ are two distinct pairs of non-adjacent vertices. This contradiction 
completes  the discussion of Case 13.1 when $l=3$.   

Assume that $l\geq 4$. Let $y_gx_l\in D$, where $g\in [1, 2]$. Then, by the minimality of $l$, the vertices $y_g, x_{l-1}$ are non-adjacent, $y_{3-g}x_{l-1}\notin D$ and $x_{l-2}y_{3-g}\notin D$. Hence by Claim 11 we get that $x_{l-2}y_g\notin D$. From the minimality of $l$ and $d^-(x_2,R)=0$ (for $l=4$) it follows that $x_{l-2}$ is not adjacent with $y_1$ and $y_2$, i.e., $d(x_{l-2},R)=0$. This together with $d^-(x_2,R)=d^-(x_{l-1},R)=0$ and Claim 3 implies that $A(R, \{x_{2}, x_{3}, \ldots , x_{l-2}\})=\emptyset$, $d^+(x_1,R)=0$, $d^-(x_1,R)\geq 1$ and $d^+(x_{l-1},R)\geq 1$. It follows that $y_2x_1\in D$ since $y_1x_1\notin D$.

Assume first that $y_g=y_2$. Then $x_{l-1}y_1\in D$. Using Lemma 3.2(ii) and 
$$
d(y_1,C[x_1,x_{l-2}])=d(y_2,C[x_2,x_{l-1}])=0$$ 
we obtain that 
$$d(y_1)=d(y_1,\{y_2\})+d(y_1,C[x_{l-1},x_{n-2}])\leq n-l+2 \quad \hbox{and}$$ 
$$ d(y_2)=d(y_2,\{y_1\})+d(y_2,C[x_{l},x_{1}])\leq n-l+2. \eqno (16)
$$
Now we extend the path $P_0:=C[x_l,x_{n-2}]$ with the vertices $x_{1}, x_{2}, \ldots , x_{l-1}$ as much as possible. Then some vertices $z_1, z_2, \ldots , z_d\in \{x_{1}, x_{2}, \ldots , x_{l-1}\}$, $d\in [2,l-1]$, are not on the extended path $P_e$. Therefore by Lemma 3.2, we have that  $d(z_i,C)\leq n+d-3$. If there is a $z_i\notin \{x_1,x_{l-1}\}$, then $d(z_i)\leq n+d-3$ and by (16), $d(z_i)+d(y_1)\leq 2n-2$ and  $d(z_i)+d(y_2)\leq 2n-2$, which contradicts Lemma 3.5 since $z_i$ is not adjacent with $y_1$ and $y_2$. Therefore assume that $\{z_1,z_2\}=\{x_1,x_{l-1}\}$ ($d=2$). Then $P_e$ ($e=l-3\geq 1$) is an $(x_l,x_{n-2})$-path with vertex set $V(C)-\{x_1,x_{l-1}\}$. Thus, we have that $y_2P_ey_1y_2$ is a cycle of length $n-2$. Therefore, by Claim 11, $x_1x_{l-1}\in D$, and hence $x_1x_{l-1}P_{e-1}y_1y_2x_1$ is a cycle of length $n-1$, which is a contradiction to our supposition.

Assume second that $y_g=y_1$. Then $x_{l-1}y_2\in D$ and $d(y_1, C[x_1,x_{l-1}])=0$. Using Lemma 3.2, we obtain that for this case (16) also holds, since $x_1y_2\notin D$ and $y_2x_{l-1}\notin D$. Again we extend the path $C[x_l,x_{n-2}]$ with  
vertices $x_{1}, x_{2}, \ldots , x_{l-1}$ as much as possible. Then some vertices $z_1, z_2, \ldots , z_d\in \{x_{1}, x_{2}, \ldots , x_{l-1}\}$, $d\in [1,l-1]$, are not on the extended path $P_e$. Similarly to the first case, we obtain that $z_i\notin \{x_2,x_3,\ldots , x_{l-2}\}$ (i.e., $z_i=x_1$ or $z_i=x_{l-1}$) and $d(z_i)\leq n+d-2$. Notice that $C':=y_1P_ey_1$ is a cycle of length $n-d-1$ with vertex set $V(C)\cup \{y_1\}-\{z_1,z_d\}$. From Claim 11 it follows that $d=2$, i.e.,     $\{z_1,z_d\}=\{x_1,x_{l-1}\}$.  From (16) and $d(z_i)\leq n+d-2$ we obtain that $d(y_1)+d(x_1)\leq 2n-2$ and $d(y_1)+d(x_{l-1})\leq 2n-2$, which contradicts Lemma  3.5, since $y_1,x_1$ and $y_1, x_{l-1}$ are two distinct pairs of non-adjacent vertices.
 \\

\textbf{Case 13.2}. $A(R\rightarrow \{x_3,x_4,\ldots , x_{n-3}\})= \emptyset$. Then 
$A(R\rightarrow \{x_{n-2},x_1\})\not=\emptyset$ since $d^-(x_2,R)=0$ and $D$ is strong, and $y_1, x_{n-3}$ are non-adjacent (by Claim 12). For this case we distinguish three subcases.

\textbf{Subcase 13.2.1}. $y_2x_{n-2}\in D$. Then it is easy to see that $d(x_{n-3},R)=0$. This together with  $y_2x_{n-2}\in D$  and Claim 3 implies that $A(\{x_1,x_2,\ldots , x_{n-3}\}\rightarrow R)=\emptyset$. Therefore  $d(R, \{x_2,x_3,\ldots , x_{n-3}\})$ $=\emptyset$ and $d(y_1)$, $d(y_2)\leq 4$ (since $y_2x_1\notin D$ by Claim 12) and $d(x_{n-3})\leq 2n-6$. From these it follows that $d(y_1)+d(x_{n-3})\leq 2n-2$ and $d(y_2)+d(x_{n-3})\leq 2n-2$, which contradicts Lemma 3.5.

\textbf{Subcase 13.2.2}. $y_2x_{n-2}\notin D$ and $y_2x_1\in D$. Then using Claim 12 it is easy to see that $y_2$ and $x_{n-2}$ are  non-adjacent.

 Let $x_{n-3}y_2\in D$. Then $y_1x_{n-2}\in D$ (by Claim 11). Using Claims 11 and 12 we 
obtain that $x_{n-4}$ is not adjacent with $y_1$ and $y_2$. From Claim 3 it follows that
 $A(\{x_1, x_2, \ldots , x_{n-4}\}\rightarrow R)=\emptyset$ and $A(R, C[x_2,x_{n-4}])=\emptyset$. Therefore and $d(y_1)=d(y_2)= 4$ and $d(x_{2})\leq 2n-6$. From these it follows that $d(y_1)+d(x_{2})\leq 2n-2$ and $d(y_2)+d(x_{2})\leq 2n-2$, which contradicts Lemma 3.5 since 
$x_2$, $y_1$ and $x_2,y_2$ are two distinct pairs of non-adjacent vertices. 

Let now $x_{n-3}y_2\notin D$. Then $y_2,x_{n-3}$ are non-adjacent and hence, $d(x_{n-3},R)=0$. Now from Claim 3 it follows that  
  $A(\{x_2, x_3, \ldots , x_{n-3}\}\rightarrow R)=\emptyset$. Therefore 
$$d(y_1,C[x_1, x_{n-3}])=d(y_2,C[x_2, x_{n-2}])=0,$$
 $d(y_1)\leq 4$, $d(y_2)\leq 4$ and $d(x_{2})\leq 2n-6$. This contradicts Lemma 3.5  since 
$x_2$, $y_1$ and $x_2,y_2$ are two distinct pairs of non-adjacent vertices. 

\textbf{Subcase 13.2.3}. $y_2x_{n-2}\notin D$ and $y_2x_1\notin D$. Then $y_1x_{n-2}\in D$ (since $D$ is strong), the vertex $y_1$ is not adjacent with vertices $x_{n-3}$, $x_{n-4}$ and
$x_{n-4}y_2\notin D$, i.e., the vertices $y_2, x_{n-4}$ also are non-adjacent. Using Claim 3, we can assume that $A(C[x_1, x_{n-4}]\rightarrow R)=\emptyset$. Therefore  $d(y_1)=4$, $d(y_2)\leq 3$ and $d(x_{1})\leq 2n-6$. This contradicts Lemma 3.5 since 
$x_1$ is not adjacent with $y_1$ and $y_2$. This completes the proof of Claim 13. \fbox \\\\

\textbf{Claim 14}. If $x_{i}y_f\in D$ and  the vertices $y_f, x_{i+1}$ are non-adjacent, then the vertices $x_{i+1}, y_{3-f}$ are adjacent, where $i\in [1,n-2]$ and $f\in [1,2]$.

\textbf{Proof}. Without loss of generality, we may assume that  $x_i=x_{n-2}$ (i.e., $x_{i+1}=x_1$) and $ y_f=y_1$. Suppose, on the contrary, that  $x_1, y_2$ are non-adjacent. From Claims 11 and 13 it follows that  $y_1x_2\notin D$ and $x_2y_1\in D$. Therefore $A(R\rightarrow \{x_1,x_2\})=\emptyset$. If $n=5$, then $x_2y_1, x_3y_1\in D$ which contradicts Claim 12. Assume therefore that $n\geq 6$. As $D$ is strong there is a vertex $x_l$ with $l\in [3, n-2]$ such that $d^-(x_l,R)\geq 1$ (say $y_gx_l\in D$) and $A(R\rightarrow C[x_1, x_{l-1}])=\emptyset$. Then the vertices $x_{l-1}, y_g$ are non-adjacent and $d(x_{l-2}, R)=0$ (by $x_{l-2}y_{3-g}\notin D$ and by Claim 11). Now, since $x_{n-2}y_1$ and $x_2y_1\in D$, there exists a vertex $x_r\in C[x_{n-2}, x_{l-3}]$ (if $l=3$, then $x_{n-2}=x_{l-3}$) such that $d^+(x_r,R)\geq 1$ and $A(R, C[x_{r+1}, x_{l-2}])=\emptyset$. This contradicts Claim 3.  Claim 14 is proved. \fbox \\\\

\textbf{Claim 15}. If $x_{i}y_j\in D$, where $i\in [1,n-2]$ and $j\in [1,2]$, then $y_jx_{i+2}\in D$. 

\textbf{Proof}.  Without loss of generality, we may assume that $x_i=x_{n-2}$ and $y_j=y_1$. 
Suppose that the claim is not true, that is $x_{n-2}y_1\in D$ and $y_1x_2\notin D$. Then, by Claims 12 and 13, the vertices  $y_1, x_1$ are non-adjacent, $x_2y_1\in D$ (hence, $n\geq 6$) and $y_1, x_3$ are also non-adjacent. From this, by Claim 14  we obtain that the vertex $y_2$ is adjacent with vertices  $x_1$ and $ x_3$. 
 Therefore either $y_2x_3\in D$ or $x_3y_2\in D$.

\textbf{Case 15.1}.  $y_2x_3\in D$. Then $x_2, y_2$ are non-adjacent (by Claim 12), $x_2x_1\in D$ and $x_1y_2\notin D$ by Claim 11 (for otherwise $D$ would has a cycle $C'$ of length $n-2$ for which $\langle V(D)-V(C') \rangle$ is not strong). Notice that $y_2x_1\in D$. Since neither $y_1$ nor $y_2$ cannot be inserted into $C$, $y_1x_2\notin D$ and $y_1, x_1$ are non-adjacent (respectively, $x_1y_2\notin D$ and $y_2, x_2$ are non-adjacent) using Lemma 3.2(ii), we obtain that
$$
d(y_1)\leq n-1 \quad \hbox{and} \quad d(y_2)\leq n-1. \eqno (17)
$$ 
Notice that $x_{n-2}x_2\notin D$ and $x_{1}x_3\notin D$. Therefore, since neither $x_1$ nor $x_2$ cannot be inserted into
$C[x_3, x_{n-2}]$ (otherwise we obtain a cycle of length $n-1$), again using Lemma 3.2(ii), we obtain that
$$
d(x_1)\leq n-1 \quad \hbox{and} \quad d(x_2)\leq n-1. \eqno (18)
$$
It is easy to check that $n\geq 7$.\\

\textbf{Remark}. Observe that from (17), (18) and Lemma 3.5 it follows that if $x_i\not=x_1$ and $y_1, x_i$ are non-adjacent or $x_i\not=x_2$ and $x_i, y_2$ are non-adjacent, then $d(x_i)\geq n$.\\

Assume first that $d^+(y_1,C [x_4,x_{n-2}])\geq 1$. Let $x_{l}$, $l\in [4,n-2]$, be the first vertex after $x_3$ that   $y_1x_{l}\in D$. Then the vertices  $y_1$ and $x_{l-1}$
are non-adjacent. Therefore $y_1$ and $x_{l-2}$ are adjacent (by Claim 13) and hence, $x_{l-2}y_1\in D$ because of $x_2y_1\in D$ and minimality of $l$ ($l-1\not=4$). Since $x_{l-1}$ cannot be inserted into $C[x_{l},x_{l-2}]$, using Lemma 3.2 and the above Remark, we obtain that $d(x_{l-1})=n$ and hence, $d(y_1)=n-1$ (by Lemma 3.5).
This together with $d(y_1, \{x_1,x_2,x_3,y_2\})=3$ implies that $d(y_1, C[x_4, x_{n-2}])=n-4$. Again using Lemma 3.2, we obtain that $y_1x_4\in D$ (since $|C[x_4,x_{n-2}]|=n-5$). Thus $y_1C[x_4,x_2]y_1$ is a cycle of length $n-2$. Therefore, $x_3y_2\in D$ (by Claim 11), $y_1x_5\notin D$ and the vertices $y_2, x_4$ are non-adjacent (by Claim 12). From $y_1x_5\notin D$ (by Lemma 3.2) we 
obtain that $d(y_1,C[x_5,x_{n-2}])\leq n-6$. Therefore $x_4 y_1\in D$ and  $d(y_1,C[x_5,x_{n-2}])= n-6$. Now it is easy to see that
$y_1, x_5$ are non-adjacent (by Claim 12) and $y_2, x_5$ are adjacent (by Claim 13). Therefore, $d(y_1,C[x_6,x_{n-2}])= n-6$ and $y_1x_6\in D$ (by Lemma 3.2), $y_2x_5, x_5y_2\in D$ (by Claim 11), $y_1x_7\notin D$ (by Claim 12). One readily sees that, by continuing the above procedure, we eventually obtain that $n$ is even and 
$$
N^-(y_1)=\{y_2,x_2,x_4,x_6,\ldots , x_{n-2}\}, \quad N^+(y_1)= \{y_2,x_4,x_6,\ldots , x_{n-2}\},
$$
$$
N^-(y_2)= \{y_1,x_3,x_5,\ldots , x_{n-3}\}, \quad N^+(y_2)=\{y_1,x_1,x_3,x_5,\ldots , x_{n-3}\}.
$$
From Claim 11 it follows that $x_ix_{i-1}\in D$ for all $i\in [4,n-2]$ and $x_2x_1\in D$. It is easy  to see that $x_1x_3\notin D$ and $x_3x_5\notin D$. Therefore, since $x_3$ cannot be inserted into $C[x_5,x_1]$, by Lemma 3.2, we have $d(x_3,C[x_5,x_{1}])\leq n-6$. This together with $d(x_3)=n$ (by Remark) implies that $d(x_3,\{x_2,x_4,y_2\})=6$. In particular, $x_3x_2\in D$. Now we consider the vertex $x_{n-2}$. Note that $d(x_{n-2})=n$ (by Remark), $x_{n-2}x_2\notin D$ and $x_{n-4}x_{n-2}\notin D$. From this it is not difficult to see that $d(x_{n-2},C[x_2,x_{n-4}])\leq n-6$ and  $x_1x_{n-2}\in D$. It follows that $x_{n-2}x_{n-3}\ldots  x_4x_3y_2x_1x_{n-2}$ is a cycle of length $n-2$, which does not contain the vertices $y_1$ and $x_2$. This contradicts Claim 11, since $y_1x_2\notin D$ (by our supposition), i.e., $\langle \{y_1,x_2\} \rangle$ is not strong.

Assume next that $d^+(y_1, \{x_4,x_5,\ldots , x_{n-2}\})=0$. Then from Claims 12 and 13 it follows that 
$$
N^-(y_1)=\{y_2,x_2,x_4,\ldots , x_{n-2}\} \quad \hbox{and} \quad N^+(y_1)=\{y_2\}. \eqno (19)
$$
By Claim 14 we have that the vertex  $y_2$ is adjacent with each vertex $x_i\in  \{x_1,x_3,\ldots , x_{n-3}\}$. It is easy to see that $x_{n-3}y_2\notin D$ and hence, $y_2x_{n-3}\in D$ (for otherwise if $x_{n-3}y_2\in D$, then $y_2C[x_1,x_{n-3}]y_2$ is a cycle of length $n-2$, but $\langle \{x_{n-2},y_1\} \rangle$ is not strong, a contradiction to Claim 11). By an argument similar to that in the proof of (19) we deduce that 
$$
N^+(y_2)=\{y_1,x_1,x_3,\ldots , x_{n-3}\} \quad \hbox{and} \quad N^-(y_2)=\{y_1\}. 
$$
Thus we have that $y_1y_2C[x_5,x_2]y_1$ is a cycle of length $n-2$ and $x_3$ cannot be inserted  into $C[x_5,x_2]$. Therefore by Lemma 3.2(ii), $d(x_3,C[x_5,x_2])\leq n-4$ since  $x_3x_5\notin D$. This together with $d(x_3,\{x_4,y_1,y_2\})\leq 3$ implies that  $d(x_3)\leq n-1$ which contradicts the above Remark  that $d(x_3)\geq n$.\\

\textbf{Case 15.2.} $y_2x_3\notin D$. Then, as noted above, $x_3y_2\in D$. Therefore $d(y_2, \{x_2, x_4\})=0$ (by Claim 12 and $y_2x_2\notin D$), $y_1x_4\notin D$ (by Claim 11), $x_4y_1\in D$ (by  Claim 14), the vertices $x_5, y_1$ are non-adjacent and the vertices $y_2, x_5 $ are adjacent (by Claim 14). Since $x_3y_2\in D$, $y_1x_4\notin D$ and $y_1, x_5$ are adjacent, from Claim 11 it follows that $y_2x_5\notin D$ and $x_5y_2\in D$. For the same reason, we deduce that  
$$
N^-(y_1)=\{y_2,x_2,x_4,\ldots , x_{n-2}\} \quad N^-(y_2)=\{y_1,x_1,x_3,\ldots , x_{n-3}\} \quad \hbox{and} \quad A(R\rightarrow V(C))=\emptyset,
$$
 which contradicts that $D$ is strong. This contradiction completes the proof of Claim 15. \fbox \\\\

We will now complete the proof of Theorem by showing that $D$ is isomorphic to $K^*_{n/2,n/2}$. 
Without loss of generality, we assume that $x_{n-2}y_1\in D$. Then using Claims 11, 12, 13   and 15 we conclude that $y_1, x_1$ are non-adjacent (Claim 12), $y_1x_2\in D$ (Claim 15), $x_1y_2, y_2x_1\in D$ (Claim 11), $x_2, y_2$ also are non-adjacent  (Claim 12) and $y_2x_3\in D$ (Claim 15). By continuing these procedure, we eventually obtain that $n$ is even and 
$$
N^+(y_1)=N^-(y_1)=\{y_2, x_2, x_4,\ldots , x_{n-2}\} \quad \hbox{and} \quad N^+(y_2)=N^-(y_2)=\{y_1, x_1, x_3,\ldots , x_{n-3}\}.
$$
If $x_ix_j\in D$ for some $x_i, x_j\in \{ x_1, x_3,\ldots , x_{n-3}\}$, then clearly $|C[x_i,x_j]|\geq 5$ and $x_ix_jx_{j+1}\ldots x_{i-1}y_1x_{i+1}$ $\ldots x_{j-2}y_2x_i$ is a cycle of length $n-1$, contrary to our assumption. Therefore $\{y_1, x_1, x_3,\ldots , x_{n-3}\}$ is a independent set of vertices. For the same reason $\{y_2, x_2, x_4,\ldots , x_{n-2}\}$ also is a independent set of vertices. Therefore $D$ is isomorphic to $K^*_{n/2,n/2}$. This completes the proof of Theorem 1.10. \fbox \\\\

\section {Concluding remarks}

A Hamiltonian bypass in a digraph is a subdigraph obtained from a Hamiltonian cycle of $D$ by reversing one arc.

Using Theorem 1.10, we have proved that if a strong digraph $D$ of order $n\geq 4$ satisfies the condition $A_0$, then $D$ contains a Hamiltonian bypass or $D$ is isomorphic to one tournament of order 5.\\


\begin{thebibliography}{25}

\bibitem{[1]} J. Bang-Jensen, G. Gutin, Digraphs: Theory,  Algorithms and Applications, Springer, 2000.

\bibitem{[2]} J. Bang-Jensen, G. Gutin, H. Li, "Sufficient conditions for a digraph to be Hamiltonian", {\it J. Graph Theory}, vol. 22 no. 2, pp. 181-187, 1996.

\bibitem{[3]} J. Bang-Jensen, Y. Guo, A.Yeo, "A new sufficient condition for a digraph to be Hamiltonian", {\it Discrete Applied Math.}, vol. 95, pp. 61-72, 1999.

\bibitem{[4]} J.A. Bondy, Basic graph theory: paths and circuits. {\it In Handbook of combinatorics}, Vol. 1, 2, Elsevier, Amsterdam, 1995.

\bibitem{[5]} J.A. Bondy, C. Thomassen, "A short proof of Meyniel's theorem", {\it Discrete Math.}, vol. 19, no. 1, pp. 195-197, 1977.
\bibitem{[6]} S.Kh. Darbinyan, "Pancyclic and panconnected digraphs" {\it Ph. D. Thesis, Institute Mathematici Akad. Nauk BSSR, Minsk, 1981} (see also, Pancyclicity of digraphs with the Meyniel condition, {\it Studia Sci. Math. Hungar.}, 20 (1-4), 95-117, (1985)) (in Russian). 

\bibitem{[7]} S.Kh. Darbinyan, "A sufficient condition for the Hamiltonian property of digraphs with  large semidegrees", {\it Akad. Nauk Armyan. SSR Dokl.}, vol. 82, no. 1, pp. 6-8, 1986 (see also arXiv: 1111.1843v1 [math.CO] 8 Nov 2011). 

\bibitem{[8]} S.Kh. Darbinyan, "On the pancyclicity of digraphs with large semidegrees", {\it Akad. Nauk Armyan. SSR Dokl.}, vol. 83, (3) (1986) 99-101 (see also arXiv: 1111.1841v1 [math.CO] 8 Nov 2011). 

\bibitem{[9]} S.Kh. Darbinyan, I.A. Karapetyan, "On longest non-Hamiltonian cycles in digraphs with the conditions of Bang-Jensen, Gutin and Li", {\it Preprint available at htte: arXiv 1207.5643v2 [math.CO]}, 20 Sep 2012. 

\bibitem{[10]} S.Kh. Darbinyan, I.A. Karapetyan, "A note on long non-Hamiltonian cycles in one class of digraphs", {\it Preprint available at htte: arXiv 1209.4456v1 [math.CO]}, 20 Sep 2012. 

\bibitem{[11]} R. H\"{a}ggkvist, C. Thomassen, "On pancyclic digraphs", {\it J. Combin. Theory Ser. B}, vol. 20(1), pp. 20-40, 1976.

\bibitem{[12]} A. Ghouila-Houri, "Une condition suffisante d'existence d'un circuit hamiltonien", {\it C. R. Acad. Sci. Paris Ser. A-B}, no. 25, pp. 495-497, 1960.

\bibitem{[13]} G. Gutin, "Characterizations of vertex pancyclic and pancyclic ordinary  complete multipartite digraphs",  {\it Discrete Math.}, vol.141, pp.153-162, 1995. 


\bibitem{[14]} Y. Manoussakis, "Directed Hamiltonian graphs", {\it J. Graph Theory}, vol. 16, no. 1, pp. 51-59, 1992. 

\bibitem{[15]} M. Meyniel, "Une condition suffisante d'existence d'un circuit hamiltonien dans un graphe oriente", {\it J. Combin. Theory Ser. B}, vol. 14, pp. 137-147, 1973.

\bibitem{[16]} C. Thomassen, An Ore-type condition implying a digraph to be pancyclic", {\it Discrete Math.}, vol.19, no 1, pp.85-92, 1977. 

\bibitem{[17]} C. Thomassen, "Long cycles in digraphs",  {\it Proc. London Math. Soc.}, vol. 3, no. 42, pp. 231-251, 1981.

\bibitem{[18]} D.R. Woodall, "Sufficient conditions for circuits in graphs", {\it Proc. London Math. Soc.}, no. 24, pp. 739-755, 1972.


\end{thebibliography}
\end{document}